\documentclass[11pt]{article}

\title{A note on comparison between Birkhoff and McShane-type
integrals for multifunctions}
\author{A. Boccuto, A. R. Sambucini}

\newtheorem{teorema}{~~Theorem}[section]
\newtheorem{lemma}[teorema]{~~Lemma}
\newtheorem{remark}[teorema]{~~Remark}

\newtheorem{corollario}[teorema]{~~Corollary}

\newtheorem{definition}[teorema]{~~Definition}

\newcommand{\enne}{\mbox{$\mathbb{N}$}}
\newcommand{\erre}{\mbox{$\mathbb{R}$}}

\newcommand{\Pd}{\mbox{$\mathcal{P}_{\Delta}$}~}

\newcommand{\fine}{~\hspace*{\fill}{$\Box$}}
\usepackage{latexsym,graphicx}
\usepackage{amsmath,amssymb}
\usepackage[latin1]{inputenc}

\begin{document}
\maketitle
\date{}

\begin{center}
Department of  Mathematics and Informatics,\\
1 via Vanvitelli,   I-06123 Perugia (ITALY)\\
e-mail:antonio.boccuto@unipg.it,
anna.sambucini@unipg.it
\end{center}

\begin{abstract}
Here we present some comparison results between Birkhoff and McShane multivalued integration.
\end{abstract}

\mbox{~} \vskip.1cm
\noindent
\normalsize

\section{Introduction}
Several notions of multivalued integral have been developed after the pioneering papers of Aumann and Debreu
in the sixties and they are used extensively in economic theory and optimal control, we cite here for example \cite{AU,uno, BY}.
In \cite{nostro} Boccuto and Sambucini and in \cite{tre,cascales2005,cascales2007,cascales2009} Cascales, Kadets and  Rodr\'{\i}guez  introduced
the McShane and  the Birkhoff multivalued integrals respectively, while in \cite{DPM, DPM2} Di Piazza and Musia\l \, introduced the Kurzweill-Henstock
one.
Since these kinds of integration lie strictly between
Bochner and Pettis integrability (both in the single-valued  and in multivalued cases) it is natural to study the possible
relationships between the Birkhoff and McShane integrals and with the other multivalued integrals,
in particular with the Pettis and Aumann Pettis studied also in \cite{BS1,BS2,EH}.
This paper is organized as follows: in Section 2 we recall some known results for the single-valued case, in Section 3 we recall  multivalued integrals, while in Section 4 we give the comparison results
between them and with respect to the Aumann integral
obtained with different types of selections.\\
The results that we obtain are the following: the McShane and the Birkhoff multivalued integrals are equivalent in Banach spaces with weak$^*$
separable dual unit ball (for example in separable Banach spaces) and in this case they agree also with the Aumann Pettis multivalued integrals
(this last comparison for the McShane  multivalued integral holds without separability assumption).

\section{The single-valued McShane and Birkhoff integrals}

 Throughout this paper   $\Omega$ is an abstract non-empty set and ${\cal T}$ is a topology on $\Omega$
 making $(\Omega,{\cal T}, \Sigma, \mu)$ a $\sigma$-finite quasi-Radon
 measure
 space which is \em outer regular, \rm namely such that
 $$\mu(B)=\inf \{ \mu(G):B \subseteq G \in {\cal T} \}\quad \mbox{for \, all  \,}
\, B \in \Sigma .$$

\noindent A \em generalized McShane partition \rm $P$ of $\Omega$
(\cite[Definitions 1A]{F2}) is a disjoint sequence $(E_i, t_i)_{i \in \enne}$
of measurable sets of finite measure, with $t_i \in \Omega$ for every $i \in \enne$ and
$\mu(\Omega \setminus \bigcup_i E_i)=0$. \\
A \em gauge \rm on $\Omega$ is a function $\Delta: \Omega \rightarrow {\cal T}$ such that $s \in \Delta(s)$ for every $s \in \Omega$.
A generalized McShane partition $(E_i,t_i)_i$ is \em  $\Delta$-fine \rm if $E_i \subset \Delta(t_i)$ for every $i \in \enne$.\\
From now on, let $X$ be a Banach space, denote by the symbol ${\cal P}$ the class of all generalized McShane partitions of $\Omega$
and by ${\cal P}_{\Delta}$  the set of all \mbox{$\Delta$}-fine elements of ${\cal P}$.

\begin{definition}\label{integral}
\rm A function $f:\Omega \rightarrow X$ is said to be
\begin{description}
\item[\ref{integral}.1] \it  McShane integrable, \rm with integral $w$,
if for every $\varepsilon >0$
there exists a gauge \mbox{$\Delta:\Omega \rightarrow {\cal T}$} such that
$$\limsup_{n \rightarrow + \infty} \, \left \| w - \sum_{i=1}^{n} \, \mu(E_i) f(t_i) \right\| \leq \varepsilon $$ for every
$\Delta$-fine McShane partition $(E_i,t_i)_i$.
In this case, we write
$\int_{\Omega} \, f =w$ (see \cite[Definition 1A]{F2});
\item[\ref{integral}.2]
\it  Birkhoff integrable, \rm  if  for every $\varepsilon > 0$ there exists
a countable partition  $\Gamma = (A_n)_n$ of $\Omega$ in $\Sigma$,
for which $f$ is summable (namely
$J(f,\Gamma) := \{ \sum_n \mu(A_n) f(t_n): t_n \in A_n \}$ is made up of unconditionally convergent series)
and $\sup_{x,y \in J(f,\Gamma)} \| x - y\| \leq \varepsilon$. In this case the Birkhoff integral of $f$\,  is
\[ (B)\int_{\Omega} f d\mu = \bigcap \{ \overline{co(J(f,\Gamma))}: f \, 
\mbox{is summable with respect to \,\,} \Gamma 
\]
(see \cite{due},
\cite[(b), Section 4]{F4}).
\end{description}
\end{definition}

The Birkhoff integral was also defined  for $\sigma$-finite measure spaces (\cite{due}) considering only
partitions into sets of finite measure. We now recall the known results on the McShane and Birkhoff
integrals in the single-valued case.

\begin{remark} \rm (See for reference \cite{F4,F2,quattro,F3,G,cinque,T84}) \label{remarkbmcs}
Let $\Omega$\, be a quasi-Radon probability space.
If $f$ is Birkhoff integrable, then $f$ is McShane
integrable and the respective integrals coincide;
 if $B_{X^*}$  is separable in the weak$^*$-topology (e.g. when $X$ is separable),
then a function $f$ is Birkhoff integrable
if and only if $f$ is McShane integrable (\cite[Theorem 10]{quattro}).
 Birkhoff integrability is in general stronger than McShane integrability,  since it is possible to construct
 a bounded McShane integrable function $f:[0,1] \to l_{\infty}([0,1])$ which is not
Birkhoff integrable (see \cite[\S 8, Example]{quattro}).
\\ If we compare these definitions of integrability with the other ones
known in the literature we remember that   
 Bochner integrability implies  McShane integrability and the two
integrals are equal (\cite[Theorem 1K]{F2}), while
McShane integrability implies  Pettis integrability  and the two integrals
coincide (\cite[Theorem 1Q]{F2}). Finally, if the Banach space $X$ is
separable, then  Birkhoff, McShane  and Pettis integrability are equivalent (\cite[Corollary 4C]{F2}, \cite{Pettis}).\\
For a more detailed investigation on the properties of the Birkhoff and Pettis integrals see also
\cite{cascales2007,cascales2005,musial}.
\end{remark}

\section{Multivalued integrals}

We skip now to the multivalued case.
Independently  in \cite{nostro} and  in \cite{tre} the authors
studied two kinds of multivalued integration related to McShane and Birkhoff integrability using
a R{\aa}dstrom embedding theorem and compared them with the usual Aumann integral.\footnote{We point out that
all the results given in \cite{nostro}  are expressed  for $\Omega=[a,b]$, where $a, b \in [-\infty, +\infty]$,
$a<b$, only for the sake of  simplicity. Moreover,  ${\cal T}$, $\Sigma$ and $\mu$ are
the families of all open subsets of $[a,b]$, the $\sigma$-algebra
of all Lebesgue measurable subsets of $[a,b]$ and the Lebesgue
measure on $[a,b]$ respectively. We observe that all the results given there   hold as well
whenever $\Omega$ is any non-empty $\sigma$-finite quasi Radon outer regular measure space.}

Let $X$ be a Banach space, $cwk(X)$~ [$ck(X)$] denote the family of all convex and weakly compact
[respectively  convex and compact]
subsets of  $X$. We denote by
the symbol $d(x,C)$ the usual distance between a point and a
nonempty set $C \subset X$, namely
$d(x,C) = \inf \{ \|x - y\|: y \in C\}$, $ \delta^* (x^*, A)= \sup_{x \in A} x^*(x)$ and with $h$ the usual Hausdorff distance.
We recall also that a multifunction $F$ is \textit{measurable} if
$F^- (C)$ is a Borel set for every closed set $C \subset X$, and that $F$ is \textit{integrably
bounded} if there exists $g \in L^1(\Omega)$ such that
\(h(F(t), \{ 0\}) \leq g(t) ~ \mu - \mbox{a.e.}~.\)\\
Thanks to the R{\aa}dstrom embedding theorem (see \cite{rad}),
$cwk(X)$ endowed with the Hausdorff distance $h$ is a complete metric space that can be isometrically embedded
 into a Banach space, for example in the Banach space of bounded real valued functions defined on $B_{X^*}$, \,
 $l_{\infty} (B_{X^*})$ endowed
with the supremum norm \mbox{$\| \cdot \|_{\infty}$} by means of the mapping
$j: cwk(X) \mapsto l_{\infty} (B_{X^*})$ given by $j(A) := \delta^* (\cdot, A)$
(see \cite[Theorems II.18 and II.19]{CV} and \cite[Lemma 1.1]{tre} for the notations).
So the authors in \cite{nostro} and \cite{tre}
 defined the  multivalued integrals  as follows:
\begin{definition}\label{hkint} \rm
Let $F:\Omega \rightarrow cwk(X)$ be a multifunction. For every $A \in \Sigma$ we say that $F$ is:
\begin{description}
\item[(\ref{hkint}.1)] \it McShane integrable \rm if   there
exists $I \in cwk(X)$ such that for every $\varepsilon
>0$ there exists a gauge $\Delta$ such that
$\limsup_n \, h(I,\sum_{i=1}^{n} F(t_i)\mu(E_i)) \leq \varepsilon $ for every
generalized \Pd McShane partition $\Pi = (E_i,t_i)_i$ of $A$.
In this case the McShane integral of $F$ on $A$ is defined by:
$I := (McS)\int_{A}  F(t) d\mu$  (\cite[Definition 1]{nostro});
\item[(\ref{hkint}.2)]
\it  Birkhoff integrable  \rm if the single-valued function
$j \circ F: \Omega \rightarrow l_{\infty} (B_{X^*})$ is Birkhoff integrable.
Since  $j(cwk(X))$ is a closed convex cone in $l_{\infty} (B_{X^*})$,
$ \int_A j \circ F d\mu \in cwk(X)$, and therefore there is a unique element $(B)\int_A F d\mu \in cwk(X)$,
 called the Birkhoff integral of $F$ on $A$, which  satisfies
$$ j\left( (B)\int_A F d\mu \right) = \int_A j \circ F d\mu$$
(\cite[Definition 2.1]{tre});
\item[(\ref{hkint}.3)] \it Aumann integrable \rm  if
\[ (A) \int_A F d\mu = \left\{ \int_A f d\mu, f \in \cal{C} \right\} \neq \emptyset,\]
where for $\mbox{$\cal{C}$}$ we  consider the following sets: $S^1_F$, $S^1_{McS}$, $S^1_{B}$, $S^1_{Pe}$
(the sets of all Bochner, McShane, Birkhoff and  Pettis integrable selections of
$F$ respectively);\\
\item[(\ref{hkint}.4)] \it Pettis integrable \rm if
$\delta^*(x^*, F): \Omega \rightarrow \erre$ given by $\delta^*(x^* ,F) (\omega) = \delta^*(x^* , F(\omega))$ is $\mu$-integrable
and if for every
$A \in \Sigma$ there exists $C_A \in \mbox{\em cwk}(X)$ such that
\[
\delta^*(x^* ,C_A) = \int_A \delta^*(x^*, F(\omega)) d\mu \mbox{ for every $x^* \in X^*$}.
\]
In this case we write
$C_A = (P) \displaystyle{\int_A} F d\mu$. See  \cite[Theorem 5.4]{EH} for a number of equivalent definitions.
\end{description}
\end{definition}

We remark moreover that,  in \cite[Corollary 2.7]{tre} it is observed that the definition of the Birkhoff integral  does not depend on the
 particular embedding used.

\section{Comparisons between multivalued integrals}
What we obtain in this paper is a comparison between the different types of
multivalued integrals introduced before.
First of all  we report some equivalence conditions for both the McShane and the Birkhoff integral:
\label{nota1} \rm \mbox{~}
\begin{description}
\item[(\ref{nota1}.a)]
If $F: \Omega \rightarrow cwk(X)$ is McShane integrable, then its integral coincides with
the $(\star)$-integral, namely
$(McS)\int_{\Omega}  F d\mu=\Phi(F,\Omega)$, where
\begin{eqnarray*}\label{int*}
\Phi(F,\Omega) &=& \{ x \in X: \forall \, \varepsilon > 0, \exists \mbox{~a gauge~} \Delta :
\mbox{~~for every generalized \Pd} \\ & &
 \mbox{ McShane partition}
(E_i,t_i)_{i\in \enne} \mbox{~there holds:~}
\\ & &
  \limsup_n \, d(x, \sum_{i=1}^{n} F(t_i) \mu(E_i)) \leq \varepsilon \}.
\end{eqnarray*}
(See \cite[Proposition 1]{nostro}). No measurability is required a priori
and so we can define the multivalued integral also in non separable Banach  spaces; moreover,
if $F$ is single-valued, then $\Phi(F, \Omega)$ coincides with the classical McShane integral, if it exists.
\item[(\ref{nota1}.b)] The Birkhoff integrability of $F$ is equivalent to:
\begin{description}
\item[{\rm (i)}] there is $W \in cwk(X)$ with the following property:
for every $\varepsilon >0$ there exists a countable partition $\Gamma_0$
of $\Omega$ in $\Sigma$ such that for every countable partition $\Gamma = \{ A_n \}$
of $\Omega$ in $\Sigma$ finer than $\Gamma_0$ and any choice of points
$t_n$ in $\Gamma$, $n \in \enne$, the series $\sum_{n=1}^{\infty} \mu(A_n) F(t_n)$
is unconditionally convergent and  $h \left( \sum_{n=1}^{\infty}
\mu(A_n) F(t_n), W \right) \leq \varepsilon.$
In this case, $\displaystyle{W=(B) \int_{\Omega} F \, d\mu}$.
\end{description}
(see \cite[Proposition 2.6]{tre}).
 Moreover in \cite[Proposition 2.9]{tre} it is showed that for bounded multifunctions
$F$ the Birkhoff integrability is equivalent to both Birkhoff and Bourgain properties.
\end{description}

As said in Remark \ref{remarkbmcs}, the Birkhoff and McShane single-valued integrals
 coincide when the Banach space has weak$^*$ separable dual unit ball. \rm
The following simple lemma is the key to extend this result from the single-valued case
to the case of  $cwk(X)$-valued functions. \\
If $A$ is a subset of a real vector space $V$,
we denote by aco$_{\mathbb{Q}}(A)$ (resp. aco$(A)$) the set of all elements
$v \in V$
that can be written as $v=\sum_{i=1}^n \, \lambda_i v_i$, $n \in \enne$,
with $v_i \in A$, $\lambda_i \in \mathbb{Q}$ (resp. $\lambda_i \in \erre$) and
$\sum_{i=1}^n |\lambda_i| \leq 1$.

\begin{lemma}\label{weaksep}
Let $Y:=C_b(B_{X^*})$  be the Banach space (with the supremum norm) of all real bounded and continuous
functions on $(B_{X^*}, \tau )$, where $\tau$ is the restriction
of the Mackey topology in $X^*$. If $B_{X^*}$ is weak$^*$ separable,
then $B_{Y^*}$ is weak$^*$ separable.
\end{lemma}
{\bf Proof:} Since $B_{X^*}$ is weak$^*$ separable, there is
a countable set $D \subset B_{X^*}$ such that ${\overline{D}}^{\tau}=
B_{X^*}$ (see the proof of \cite[Lemma 3.6]{tre}). Given $d \in D$, let us consider the
element $y^*_d \in B_{Y^*}$ defined by $y^*_d(f):=f(d)$. Since $\{ y_d^*: d \in D \}
\subset B_{Y^*}$ is norming, by applying the Hahn-Banach theorem we get
$$\overline{{\rm aco}(\{y^*_d:d \in D\})}^{{\rm weak}^*}=B_{Y^*},$$ and hence
$$\overline{{\rm aco}_{\mathbb{Q}}(\{y^*_d:d \in D\})}^{{\rm weak}^*}=B_{Y^*}.$$
Therefore ${\rm aco}_{\mathbb{Q}}(\{y^*_d:d \in D\})$ is a countable weak$^*$
dense subset of $B_{Y^*}$, and the proof is complete. $\quad \Box$

\begin{corollario}\label{confronto}
Let $\Omega$ be a quasi-Radon probability space and let
$X$ be a Banach space such that $B_{X^*}$ is weak$^*$ separable. Then a multi-valued
function $F:\Omega \to cwk(X)$ is McShane integrable if and only if
$F$ is Birkhoff integrable. In this case, the two integrals coincide.
\end{corollario}
{\bf Proof:} Let $j: cwk(X) \to l_{\infty}(B_{X^*})$ be the embedding
used in \cite{tre}. Then it is easy to see that
\mbox{$F:\Omega \to cwk(X)$} is McShane integrable if and only if the
single-valued function $j \circ F: \Omega \to l_{\infty}(B_{X^*})$ is
McShane integrable according to \cite[Definitions 1A]{F2}, and in this case 
\mbox{$j\left(\int F\right)=\int j \circ F$}. Since $j(cwk(X)) \subset
C_b(B_{X^*})$, the function $j \circ F$ takes values in $C_b(B_{X^*})$,
which is a closed subspace of $l_{\infty}(B_{X^*})$ with weak$^*$ separable
dual unit ball by virtue of Lemma \ref{weaksep}.
By \cite[Theorem 10]{quattro}  $j \circ F$ is McShane integrable if and only if $j \circ F$
is Birkhoff integrable and the respective integrals coincide.  $\quad \Box$ \vspace{3mm}

The same conclusion  of Corollary \ref{confronto} can be obtained for $\sigma$-finite quasi-Radon outer
regular measure spaces.

We now want  to compare  the McShane and the Birkhoff multivalued integrals with the Aumann integral, when the
multifunction $F$ has some kind of measurability. For the case of Birkhoff integrability the result is given in
\cite[Proposition 3.1]{tre}, for the McShane case we have:

\rm

\begin{teorema}\label{treuno}
Let $F:\Omega \to cwk(X)$ be a McShane integrable multifunction,
then $F$ is Pettis integrable and for every $A \in \Sigma$ we have
\begin{eqnarray}\label{a1}
(McS) \int_A F d\mu = \overline{
\left\{ \int_A f d\mu, \,\,\, f \in S^1_{Pe} \right\}.
}
\end{eqnarray}
Moreover, if  $(\Omega, {\cal T}, \Sigma, \mu)$ is a Radon measure space 
or there is no real-valued-measurable cardinal
and every Pettis integrable selection  $f$ is measurable (that is, $f^{-1}(G) \in \Sigma$
for every norm-open set $G \subseteq X$), then
\begin{eqnarray}\label{a2}
(McS) \int_A F d\mu = \overline{ \left\{ \int_A f d\mu, \,\,\, f \in S^1_{McS} \right\} }.
\end{eqnarray}
If $X$ is separable then
\begin{eqnarray}\label{a3}
(B) \int_A F d\mu &=&  \left\{ \int_A f d\mu, \,\,\, f \in S^1_{B} \right\}  = \left\{ \int_A f d\mu, \,\,\, f \in S^1_{McS} \right\} \\ &=&
 (McS) \int_A F d\mu. \nonumber
\end{eqnarray}
Finally, if $F$ is  measurable and integrably bounded  and $X$ is separable and there exists a countable family
$(x^{*}_n)_n$ in $X^{*}$ which separates points of $X$, then
\begin{eqnarray}\label{a4}
(McS) \int_A F d\mu = \left\{ \int_A f d\mu, \,\,\, f \in S^1_{F} \right\}.
\end{eqnarray}
\end{teorema}
{\bf Proof:}
If $F$ is McShane integrable then it means that $j \circ F$ is McShane integrable thanks to the R{\aa}dstrom embedding theorem.
So  $j \circ F$ is Pettis integrable,
and then, by \cite[Proposition 4.4]{cascales2009}, $F$ is Pettis integrable and
for all $A \in \Sigma$ and  $x^*  \in B_{X^*}$ we get:
\begin{eqnarray*}
\delta^*(x^*,F)&=&\langle e_{x^*}, j \circ F \rangle \in L^1(\mu);
\end{eqnarray*}
where $e_{x^*} \in B_{l_{\infty} (B_{X^*})^*}$ is defined by: $\langle e_{x^*}, g\rangle := g(x^*)$ for every $g \in l_{\infty} (B_{X^*})$. Then
\begin{eqnarray*}
\delta^* \left( x^*,(McS) \int_A F \, d\mu \right) &=&
\langle e_{x^*}, j \left( (McS) \int_A F \, d\mu  \right) \rangle
= \langle e_{x^*}, \int_A j \circ F \, d\mu \rangle \\
&=& \int_A \langle e_{x^*},
j \circ F \rangle \, d\mu = \int_A \delta^* (x^*,F) \, d\mu.
\end{eqnarray*}
Thus for all $A \in \Sigma$ we have
$$(P) \int_A F \, d\mu=  (McS) \int_A F \, d\mu.$$
Now by \cite[Theorems 2.5 and 2.6]{cascales2009}
$F$ admits Pettis integrable selections and
\begin{eqnarray*}
(P) \int_A F \, d\mu= \overline{ \left\{ \int_A f d\mu, \,\,\, f \in S^1_{Pe} \right\} },
\end{eqnarray*}
which proves (\ref{a1}).
Now, by virtue of \cite[Theorem 1Q and Corollary 4D]{F2},
in our context McShane and Pettis integrability coincide for single-valued
functions and this proves (\ref{a2}).\\
Observe that the concept of no real-valued measurable cardinal which appears in \cite[Corollary 4D]{F2} and in
\cite[(e)]{F0} is contained in \cite[\S 438]{Fvol4} using the new terminology measure-free cardinals.
In this case all metric spaces are Radon and $f$ is McShane integrable using \cite[438D]{Fvol4} 
and \cite[Corollary 2G]{F2}.
\\
If $X$ is separable, the first equality in (\ref{a3}) for the Birkhoff integral is given in \cite[Proposition 3.1]{tre}
and the last equalities follow in an analogous way and taking into account the equivalence among Pettis, Birkhoff and McShane integrability.
So we obtain again the equivalence between the Birkhoff and the McShane multivalued integrals  in a different way.
Finally (\ref{a4}) is given in \cite[Theorem 1]{nostro}.\fine \\

In the end the comparison with the Debreu integral is obvious thanks to the given definitions
(see \cite[Proposition 3.1 (i), Theorem 3.2]{tre} and \cite[page 321 and Corollary 1]{nostro}).
Comparisons between Aumann and Debreu integrals are given also in \cite{jmaa,zaa,Sambucini,S1}.

\medskip
\noindent {\bf Acknowledgment}. The authors would like to thank D. H. Fremlin and 
J. Rodr\'{\i}guez for their helpful discussions and  suggestions given during the writing of the  paper.


\begin{thebibliography}{99}
\bibitem{AU} R. J. Aumann, {\it Integrals of set-valued functions}, 
J. Math. Anal. Appl. {\bf 12}  (1965), 1-12.

\bibitem{uno} Z. Artstein and J. A. Burns, {\it Integration of compact set-valued
functions},  Pacific J. Math. {\bf 58}  (1975), 297-307.

\bibitem{BS1}
E. J. Balder and A.R. Sambucini  {\it A note on strong convergence for Pettis integrable function},
 Vietnam J.  Math.  {\bf 31},  N. 3 (2003), 341-347.

\bibitem{BS2}
E. J. Balder and A.R. Sambucini {\it On weak compactness and lower closure results for Pettis integrable
(multi)functions},  Bull. Pol. Acad. Sci. Math.  {\bf 52},  N. 1 (2004), 53-61.

\bibitem{nostro} A. Boccuto and A. R. Sambucini,
{ \it A McShane integral for multifunctions}, J. Concr.  Appl. Math.  {\bf  2},
\rm  N. 4 (2004), 307-325.

\bibitem{BY} C. L. Byrne, {\it Remarks on the Set-Valued Integrals of
Debreu and Aumannn},  J. Math. Anal. Appl. {\bf 62} (1978), 243-246.

\bibitem{due} G. Birkhoff, {\it Integration of functions with values in a
Banach space},  Trans. Amer. Math. Soc. {\bf 38}  (1935), 357-378.

\bibitem{cascales2007} B. Cascales, V. Kadets and J. Rodr\'{\i}guez, {\it The Pettis integral for
multi-valued functions via single-valued functions},
 J. Math. Anal. Appl.  {\bf 332} (2007), 1-10.

\bibitem{cascales2009} B. Cascales, V. Kadets and J. Rodr\'{\i}guez, {\it Measurable
selectors and set-values Pettis integral in non-separable Banach spaces},
 J. Funct. Anal.  {\bf 256} (2009), 673-699.

\bibitem{tre} B. Cascales and J. Rodr\'{\i}guez,  {\it Birkhoff integral for
multi-valued functions},  J. Math. Anal. Appl. {\bf 297} (2004), 540-560.

\bibitem{cascales2005} B. Cascales and J. Rodr\'{\i}guez, {\it The Birkhoff integral and
the property of Bourgain},  Math. Ann.  {\bf 331} (2005), 259-279.

\bibitem{CV} C. Castaing and M. Valadier, {\em Convex Analysis and
Measurable Multifunctions},  Lecture Notes in Math. {\bf 580},
 Springer-Verlag  (1977).

\bibitem{DPM} L. Di Piazza and K. Musia{\l}, {\it Set-Valued Kurzweil-Henstock-Pettis Integral},
Set-Valued Anal.  {\bf 13},  (2005), 167-179.

\bibitem{DPM2} L. Di Piazza and K. Musia{\l}, {\it A Decomposition Theorem for Compact-Valued
Henstock Integral},  Monatsh. Math.  {\bf 148},  (2006), 119-126.

\bibitem{EH} K. El Amri and C.~Hess, {\it On the Pettis integral of closed
valued multifunctions},  Set-Valued Anal. {\bf 8} (2000), 329-360.

\bibitem{F0} D. H. Fremlin,
{\it Measurable functions and almost continuos functions}, Manuscripta Math. {\bf 33}
(1981), 387-405.

\bibitem{F4} D. H. Fremlin,
{\it Integration of vector-valued functions},  Atti Semin. Mat. Fis. Univ. Modena
{\bf 42} (1994), 205-211.

\bibitem{F2} D. H. Fremlin, 
{\it The generalized McShane integral},
 Illinois J. Math. {\bf 39} (1995), 39-67.

\bibitem{quattro} D. H. Fremlin, {\it The McShane and Birkhoff integrals
of vector-valued functions}, University of Essex, Mathematics Department
Research Report 92-10, version 18.5.2007, available at URL
http://www.essex.ac.uk/maths/people/fremlin/preprints.htm

\bibitem{Fvol4} D. H. Fremlin, 
{\em Measure Theory volume four: topological measure spaces}
Torres Fremlin, 25 Ireton Road, Colchester CO33AT, England (2003).

\bibitem{F3} D. H. Fremlin and J. Mendoza,
{\it On the integration of vector-valued functions},  Illinois J.
Math. {\bf 38}  (1994), 127-147.

\bibitem{G} R. Gordon, {\it The McShane Integral of Banach valued functions},
 Illinois J. Math. {\bf 34} (1990), 557-567.


\bibitem{jmaa} A. Martellotti and A. R. Sambucini,  {\it On the comparison
between Aumann and Bochner integrals},  J. Math. Anal.
Appl. {\bf 260}  N. 1 (2001), 6-17.

\bibitem{zaa} A. Martellotti and A. R. Sambucini,  {\it The finitely additive integral
of multifunctions with  closed and convex values},
Z. Anal  Anwend.   {\bf  21} N. 4 (2002), 851-864.

\bibitem{cinque} E. J. McShane, {\it A Riemann-type integral that
includes Lebesgue-Stieltjes, Bochner and stochastic integrals}.
Mem. Amer. Math. Soc. {\bf 88}, Amer. Math. Soc., Providence, R. I. (1969).

\bibitem{musial} K. Musia\l, \,
{\it Topics in the theory of Pettis integration}, \rm Rend. Istit. Mat. Univ. Trieste,
\rm \textbf{23}, \rm (1991), 177-262.

\bibitem{Pettis} B.J. Pettis, 
{\it On integration in vector spaces}.  Trans. Amer. Math. Soc., {\bf 44}, \rm (1938),
n. 2, 277-304.

\bibitem{rad} H. R{\aa}dstrom, {\it An Embedding Theorem for Spaces of Convex
Sets},  Proc. Amer. Math. Soc. {\bf 3} (1952), 165-169.

\bibitem{Sambucini} A. R. Sambucini,  {\it Remarks on set valued integrals of
multifunctions with non empty, bounded, closed and convex values}, 
Comment. Math. Prace Mat.  {\bf 39} (1999), 153-165.

\bibitem{S1} A. R. Sambucini, {\it A survey on multivalued integration},  Atti Semin. Mat. Fis. Univ. Modena 
{\bf 50} (2002), 53-63.

\bibitem{T84} M. Talagrand, {\em Pettis integral and measure theory}.
 Mem. Amer. Math. Soc. {\bf 307}, Am. Math. Soc., Providence, R. I. (1984).
 \end{thebibliography}
\end{document}